\documentclass[preprint,12pt]{elsarticle}

\usepackage{amssymb}
\usepackage{amsmath}
\usepackage{amsthm}

\usepackage{mathrsfs} 
\usepackage{subfigure}
\usepackage{subfiles}

\newtheorem{theorem}{Theorem}[section]
\newtheorem{lemma}{Lemma}[section]

\newtheorem{assumption}{Assumption}[section]

\theoremstyle{definition}

\newtheorem{remark}{Remark}[section]
\def\dfrac{\displaystyle\frac}

\journal{Systems and Control Letters}

\begin{document}

\begin{frontmatter}



\title{Output regulation for an unstable wave equation with output delay and one measurement only}


\author[s1]{Shen Wang}

\author[s1]{Zhong-Jie Han\corref{1}}
\cortext[1]{Corresponding author}
\ead{zjhan@tju.edu.cn} 

\affiliation[s1]{
organization={ School of Mathematics and KL-AAGDM, Tianjin University},
city={Tianjin},
postcode={300350}, 
country={China}}

\author[s2]{ Shuangxi Huang}

\affiliation[s2]{
organization={School of Mathematics and Statistics, Shandong Normal University}, 
city={Jinan},
postcode={250014}, 
country={China}
}

\author[s3,s4]{Zhi-Xue Zhao}

\affiliation[s3]{
organization={School of Mathematical Sciences, Tianjin Normal University},
city={Tianjin},
postcode={300387}, 
country={China}
}

\affiliation[s4]{
organization={
Institute of Mathematics and Interdisciplinary Sciences, Tianjin Normal University},
city={Tianjin},
postcode={300387}, 
country={China}
}

\begin{abstract}
This paper addresses the output regulation problem for a one-dimensional unstable wave equation subject to output delay and all-channel disturbances with unknown frequencies and amplitudes.
First, this problem is transformed into a stabilization problem for an unstable wave equation with output delay and disturbances by employing regulator equations. 
Subsequently, a backstepping-based feedforward regulator is proposed to exponentially stabilize this system.
To track the states of the unstable wave equation, the time interval is partitioned into two segments.  The observers and predictors are designed at these distinct intervals, respectively.
Therein, the observers comprise two components: a state observer proposed via dynamic compensators and an adaptive observer designed by the adaptive internal model method.
Finally, a novel error-based feedback controller is derived using a single measurement, ensuring exponential convergence of the tracking error to zero.
This work establishes the pioneering solution to the output regulation problem for distributed parameter systems (DPS) with output delay.
Numerical simulations are provided to illustrate the results.
\end{abstract}


\end{frontmatter}

\section{Introduction}

\subsection{Background}
Time delay is inherently present in practical engineering systems spanning aerospace, chemical processes, and laser models \cite{9,10}. 
It typically results in suboptimal performance and instability in distributed parameter systems (DPS), thereby posing fundamental challenges to control design \cite{11,12}.
Over the past few years, control strategies for time-delayed DPS have garnered sustained interest, yielding numerous innovative findings (see \cite{13,14,15,16,27}). 
As a prime form of time delay, output delay exists between the controller to be implemented and the information obtained through system observation. 
Existing solutions predominantly focus on canceling output delay impacts through observer-predictor architectures. 
For instance, stabilization problems for beam, Schrödinger, and wave equations have been addressed in  \cite{17,18,19}, where observers and predictors are separately proposed during the measurement-available and unavailable time intervals.
Regarding the abstract equations framework for DPS with output delay, 
In \cite{21,22}, Guo and Mei addressed the stabilization problems for systems featuring bounded control or observation operators, as well as for well-posed and regular linear systems, respectively.
In addition, Yang \emph{et al.} employed the backstepping approach to exponentially stabilize an unstable wave equation with output delay by two measurements (see \cite{8}).
Despite these advances, the output regulation problem for an unstable wave equation with output delay remains unexplored, so this work bridges this gap by establishing an innovative solution. 

It is widely recognized that output regulation problems can be transformed into stabilization problems (see \cite{23}). 
Within the domain of DPS, substantial progress has been achieved in stabilization problems, as evidenced by recent advances in \cite{31,32,33,34}. 
Concerning the stabilization of an unstable wave equation, specifically, two predominant frameworks have emerged: the backstepping approach and the dynamic compensator. 
Utilizing the backstepping approach, Krstic \emph{et al.} proposed a controller that employs a Volterra transformation to map the unstable wave equation into an exponentially stable target system (see \cite{1}). 
Additionally, infinite-dimensional dynamic compensators were embedded as auxiliary subsystems to exponentially stabilize unstable wave equations in \cite{2,7}.

Output regulation has emerged as a pivotal research field of control theory with wide-ranging applicability across engineering domains \cite{25}.
Specifically, its core objective lies in designing controllers for uncertain dynamical systems that simultaneously ensure internal stability and asymptotic output tracking of reference signals despite disturbances \cite{24}.
Significant progress has been achieved in extending output regulation to DPS such as wave equations.
Regarding conservative open-loop wave equations, both situations wherein the frequency of disturbances are known and unknown have been considered (see \cite{28,26,3}).
Notably, the exponential convergence of the tracking error in \cite{3} necessitates two measurements (tracking error and its derivative).
Regarding unstable wave equations, Guo \emph{et al.} established asymptotic stability via a backstepping-based adaptive regulator requiring two measurements (see \cite{29}),
and Li \emph{et al.} developed an effective feedback controller to achieve asymptotic convergence of the tracking error via the backstepping approach and the internal model principle (see \cite{35}). 
Despite these advances, output regulation for unstable wave equations with output delay and unknown disturbance frequencies remains an open challenge—a critical gap this work aims to resolve.

\subsection{Contributions and novelties}
This paper resolves the output regulation for an unstable wave equation with output delay and unknown frequency.
By dividing the time intervals and performing observers and predictors on distinct intervals, we cancel the negative impact of the output delay.
Based on the adaptive internal model method, we propose an adaptive observer to track the unknown frequency.
Utilizing the backstepping approach and a dynamic compensator, we derive a controller that enables the tracking error to exponentially converge to zero using a single measurement.
This work constitutes the first systematic solution bridging output delay, adaptive observer, and the backstepping approach for output regulation of time-delayed DPS.
Our contribution, compared with the existing works, can be summarized as follows:

1) This work presents the first systematic study on output regulation for an unstable wave equation with output delay.
Existing studies have primarily focused on output regulation for time-delayed systems (see \cite{37,38} and references therein), 
particularly those addressing output delays for lumped parameter systems (LPS) (see \cite{36,39}), 
while research on time-delayed DPS has predominantly examined input delays (e.g., \cite{30}). 
However, the analogous problem for output-delayed DPS has received considerably less attention.
To eliminate the effect of output delay, we design observers and predictors on different intervals, enabling exponential tracking states of the unstable wave equation.

2) This work achieves exponential convergence of the tracking error via a single measurement—a critical advancement over existing multi-measurement results. 
For the unstable wave equation, our approach not only reduces the number of measurements but also elevates the stability performance through the novel observer architecture synthesized with an infinite-dimensional dynamic compensator.
Concurrently, this dynamic compensator introduces inherent difficulties in the controller design process.

3) This work addresses the challenge where the amplitudes and frequency of disturbances are all unknown—a situation that classical observers (e.g., Luenberger observers) struggle to address.
To solve this problem, we develop an adaptive observer capable of tracking the unknown frequency in real-time, utilizing the adaptive internal model method.

\subsection{Organization}
This paper is organized as follows.
Section 2 is divided into four parts: the problem statement, feedforward regulator design, error-based observers and predictors design, and error-based feedback control design.
Section 3 gives proofs of lemmas and theorems. 
Simulation results and conclusions are provided in Section 4 and Section 5, respectively. 

\section{Main results}

\subsection{Problem statement}
In this paper, we consider a one-dimensional unstable wave equation with output delay
\begin{equation}\label{1}
\left
\{
\begin{array}{ll}
w_{t t}(x, t)=w_{x x}(x, t)+p_1(x)v(t), &
0<x<1, t>0, \\
w_x(0, t)=-q w(0, t)+p_2v(t), &
t\geq 0, \\
w(1, t)=u(t)+p_3v(t), &
t\geq 0,  \\
w(x, 0)=w_0(x),
\; w_t(x, 0)=w_1(x),&
0\leq x\leq 1, \\
y_p(t)=w(0, t-\tau), &
t\geq 0,
\end{array}
\right.
\end{equation}
where $u(t)$ is the control input and $y_p(t)$ is the output to be regulated.
We assume that $y_p(t)$ satisfies $y_p(t)\equiv0,\;0<t<\tau$, where $\tau$ is the constant time delay.
The destabilizing boundary feedback term with $q>0$ at $x=0$ acts like a spring with a negative spring constant. 
In addition, $p_1(\cdot)\in L^\infty((0,1);\mathbb{R}^{1\times 2})$ and  
$p_j \in \mathbb{R}^{1\times 2}, j=2,3$ are all unknown coefficients of the external disturbance $v(t)$.
We consider \eqref{1} in the state space $\mathcal{H}_1 = H^1_R(0,1)\times L^2(0,1)$, where 
$H^1_R(0,1)=\{ w\in H^1(0,1) | w(1)=0\}$.
The disturbance $v(t)$ and the reference signal $y_{\rm{ref}}(t)$ are both generated by the following exosystem:
\begin{equation}\label{2}
\left
\{
\begin{array}{ll}
\dot{v}(t)=Sv(t), &
t > 0, \\
v(0)=v_0, \\
y_{\rm{ref}}(t)=p_4v(t), & 
t\geq 0,
\end{array}
\right.
\end{equation}
where $S\in\mathbb{R}^{2\times 2}$,  $v_0\in\mathbb{R}^{2\times1}$, $p_4\in\mathbb{R}^{1\times 2}$ and $y_{\rm{ref}}(t)$ are all unknown. 
\begin{assumption}\label{as1}
The eigenvalues of the matrix $S$ are $\pm i\omega$, where $\omega$ is an unknown positive parameter.
\end{assumption}
We consider the coupled system \eqref{1} and \eqref{2} in the energy state space 
$\mathcal{H}_2 = \mathcal{H}_1 \times \mathbb{R}^{2}$. 
Although the output $y_p(t)$ and the reference signal $y_{\rm{ref}}(t)$ are both unknown, the tracking error 
\begin{equation}\label{3}
e(t):= y_p(t)-y_{\rm{ref}}(t) = w(0,t-\tau)-p_4 v(t) ,
\end{equation}
is assumed to be known.
The objective of the output regulation problem for \eqref{1} and \eqref{2} is to propose a tracking error-based feedback controller $u(t)$, so that $e(t)$ exponentially converges to $0$, i.e., 
\begin{equation}\label{4}
\lim_{t\to\infty}|e(t)|  
= \lim_{t\to\infty} |y_p(t)-y_{\rm{ref}}(t)|
= 0, 
\; {\rm{exponentially}}.
\end{equation}

\subsection{Feedforward regulator design}
In this subsection, we first transform the output regulation problem for \eqref{1} and \eqref{2} into a stabilization problem by regulator equations.
Then, we design a feedforward regulator for this stabilization problem using the backstepping approach.

We propose the following transformation for \eqref{1}:
\begin{equation}\label{5}
\varepsilon(x,t) = w(x,t) - \Pi(x)v(t),
\end{equation}
where $\Pi(x)$ satisfies the regulator equations
\begin{equation}\label{6}
\left\{
\begin{array}{l}
\Pi^{\prime\prime}(x)=\Pi(x) S^2 - p_1(x),\\
\Pi^\prime(0)=-q \Pi(0) + p_2,\\
\Pi (0)=p_4 e^{\tau S}.
\end{array}
\right.
\end{equation}
\begin{lemma}\label{le1}
There exists a unique solution $\Pi \in C^2( (0,1);\mathbb{R}^{1 \times 2})$ to the initial value problem \eqref{6}. 
\end{lemma}
Thus, together with \eqref{5} and \eqref{6}, we transform \eqref{1} into 
\begin{equation}\label{9}
\left
\{
\begin{array}{l}
\varepsilon_{tt}(x,t)
=\varepsilon_{x x}(x, t),
\; t>0,\\
\varepsilon_x(0, t)
=-q \varepsilon(0, t),\\
\varepsilon(1, t)
=u(t)
-\gamma_1 v(t),
\end{array}
\right.
\end{equation}
where $ \gamma_1=\Pi (1)-p_3$ and the tracking error
\begin{equation}\label{88}
e(t)=\varepsilon(0,t-\tau), \; t>\tau.
\end{equation}
The system \eqref{9} is considered in the state space $\mathcal{H}_3=H^1(0,1)\times L^2(0,1)$.
Next, we design a feedforward regulator $U(t)$ to stabilize \eqref{9} by the backstepping approach.
According to \cite{1}, we propose a transformation
\begin{equation}\label{10}
\bar{\varepsilon}(x,t)
=[(I+\mathbb{P}) \varepsilon] (x,t)
=\varepsilon(x,t) + (c_0+q)\int_0^x e^{q(x-h)} \varepsilon(h,t) dh, 
\end{equation}
where $\mathbb{P}$ is a Volterra transformation, and $c_0>0$ is a design parameter. 
It should be noted that this transformation is bounded, and invertible;
and the inverse of $(I+\mathbb{P})$ is 
\begin{equation}\label{11}
\varepsilon(x,t)
=[(I+\mathbb{P})^{-1} \bar{\varepsilon}] (x,t)
=\bar{\varepsilon}(x,t) - (c_0+q) \int_0^x e^{-c_0 (x-h)} \bar{\varepsilon}(h,t) dh.
\end{equation}
Under the transformation \eqref{10}, the system \eqref{9} is transformed into the following target system:
\begin{equation}\label{12}
\left
\{
\begin{array}{l}
\bar{\varepsilon}_{tt}(x,t)
=\bar{\varepsilon}_{x x}(x, t), \\
\bar{\varepsilon}_x(0, t)
=c_0 \bar{\varepsilon}(0,t),\\
\bar{\varepsilon}_x(1, t)
=-c_1 \bar{\varepsilon}_t (1,t),
\end{array}
\right.
\end{equation}
where $c_1>0$ is a design parameter, and the controller $u(t)$ is given in
\begin{equation}\label{13}
u(t)=
-\dfrac{1}{c_0+q} [ \varepsilon_x(1,t) + c_1 \varepsilon_t(1,t) ]
-\int_0^1 e^{q(1-h)} [ q \varepsilon(h,t) + c_1\varepsilon_t(h,t) ]dh 
+\gamma_1 v(t).
\end{equation}
Similar to the proof of \cite{1}, we can obtain the following result.
\begin{lemma}\label{le2}
Suppose that $c_0,c_1>0$ and the controller $U(t)$ is given in \eqref{13}.
For any initial value $(\varepsilon(\cdot,0), \varepsilon_t(\cdot,0)) \in \mathcal{H}_3$,
there exists a unique solution $(\varepsilon,\varepsilon_t) \in C([0,\infty);\mathcal{H}_3)$ to \eqref{9} such that 
\begin{equation}\label{14}
\| ( \varepsilon(\cdot,t), \varepsilon_t(\cdot,t)) \|_{\mathcal{H}_3} 
\leq M_1 e^{-\mu_1 t}
\|( \varepsilon(\cdot,0), \varepsilon_t(\cdot,0) )\|_{\mathcal{H}_3},
\end{equation}
for some $M_1,\mu_1>0$. 
\end{lemma}
According to Assumption \ref{as1}, the last term $\gamma_1 v(t)$ of \eqref{13} can be written as
\begin{equation}\label{15}
\gamma_1 v(t) = A\cos\omega t + B\sin\omega t,
\end{equation}
where $A,B,\omega$ are unknown parameters and satisfy $A^2+B^2>0$.
By introducing a new variable
\begin{equation}\label{16}
\eta(t) = [A\cos\omega t+B\sin\omega t \quad B\cos\omega t-A\sin\omega t]^\top,
\end{equation}
we rewrite $\gamma_1 v(t)$ as
\begin{equation}\label{17}
\left\{
\begin{array}{l}
\dot{\eta}(t)=S_\eta \eta(t),\\
\gamma_1 v(t)=\gamma_\eta \eta(t),
\end{array}
\right.
\end{equation}
where
\begin{equation}\label{18}
S_\eta
=\left[
\begin{array}{cc}
0 & \omega \\
-\omega & 0
\end{array}
\right], \;
\gamma_\eta = [1 \quad 0],\;
\eta(0) = [A \quad B]^{\top} .
\end{equation}
By \eqref{17}$_2$, the feedforward regulator \eqref{13} can be rewritten as
\begin{equation}\label{19}
u(t)
= -\dfrac{1}{c_0+q} [ \varepsilon_x(1,t) + c_1 \varepsilon_t(1,t) ]
-\int_0^1 e^{q(1-h)} [ q \varepsilon(h,t) + c_1\varepsilon_t(h,t)dh]
+\gamma_\eta \eta(t).
\end{equation}

\subsection{Error-based observers and predictors design}
In this subsection, we design observers and a predictor to recover the states of the coupled system \eqref{9} and \eqref{17}.
For any given $t>\tau$, we divide the time span $[0,t]$ into two parts: $[0,t-\tau]$ and $[t-\tau,t]$.
When $s\in[0,t-\tau]$ and $s\in[t-\tau,t]$, the states of \eqref{9} and \eqref{17} are estimated by the observers and the predictor, respectively.

\subsubsection{Observers design}
In this part, we design observers to estimate \eqref{9} and \eqref{17} by the available observation $e(s+\tau) = \varepsilon(0,s)$ for $s\in[0,t-\tau], \; t>\tau$.
First, we introduce a dynamic compensator
\begin{equation}\label{20}
\left\{
\begin{array}{l} 
Y_{1s}(x,s)=-Y_{1x}(x,s), \; 0<s<t-\tau, \\
Y_1(0,s)=-c_2e(s+\tau), \\ 
Y_1(x,0)=Y_{10}(x),
\end{array}
\right.
\end{equation}
where $0<c_2<1$ is a design parameter, and \eqref{20} is completely determined by the tracking error $e(s+\tau)$.

It should be noted that it is difficult to directly propose an observer for the PDE-ODE coupled system \eqref{9} and \eqref{17} only based on the measurement $e(s+\tau), \; s\in[0,t-\tau], t>\tau$.
In the following, we first decouple \eqref{9} and \eqref{17} by invertible transformations to address this problem.
Then, we design a state observer and an adaptive observer for the PDE and ODE parts, respectively.
The dynamic compensator \eqref{20} assists in the design of the state observer for the PDE part.

We propose the following transformations:
\begin{equation}\label{21}
\left\{
\begin{array}{l}
z(x,s)=\varepsilon(x,s)-g_1(x)\eta(s),\\
Y_2(x,s)=Y_1(x,s)-g_2(x)\eta(s), 
\end{array}
\right.
\end{equation}
where $g_1(x)$ and $g_2(x)$ satisfy
\begin{equation}\label{22}
\left\{
\begin{array}{l} 
g_1^{\prime\prime}(x)
=g_1(x)S_\eta^2,\\
g_1^{\prime}(0)=0,\\
g_1(1)= -\gamma_\eta-g_2(1),\\
g_2^{\prime}(x)
=-g_2(x)S_\eta,\\
g_2(0)=-c_2g_1(0).
\end{array}
\right.
\end{equation}
\begin{lemma}\label{le3}
There exists a unique solution
$g_1 \in C^2([0,1];\mathbb{R}^{1\times 2})$ and 
$g_2 \in C^1([0,1]; $ $ \mathbb{R}^{1\times 2})$ to the system \eqref{22}.
\end{lemma}

Together with \eqref{9}, \eqref{17}, \eqref{20}, and \eqref{21}, the extended $(z,Y_2)-$system is governed by
\begin{equation}\label{27}
\left\{
\begin{array}{l}
z_{s s}(x, s) = z_{x x}(x, s), \; 0<s<t-\tau,\\
z_x(0, s) = -qe(s+\tau), \\
z(1, s) = u(s)-(Y_2(1,s)-Y_1(1,s)), \\
Y_{2s}(x, s) = -Y_{2x}(x,s), \; 0<s<t-\tau,\\ 
Y_2(0,s) = -c_2 z(0,s), \\ 
e(s+\tau) = z(0,s)+g_1(0)\eta(s).
\end{array}
\right.
\end{equation}
A remarkable feature of \eqref{27} is that the disturbance $\eta(s)$ appears in tracking error $e(s+\tau)$ only, and 
thereby we have decoupled the systems \eqref{9} and \eqref{17}.
Next, we design the following observer for \eqref{27}:
\begin{equation}\label{28}
\left\{
\begin{array}{l}
\hat{z}_{s s}(x, s) = \hat{z}_{x x}(x, s), \\
\hat{z}_x(0, s)=-qe(s+\tau), \\
\hat{z}(1, s) = u(s)-(\hat{Y}_2(1,s)-Y_1(1,s)), \\
\hat{Y}_{2s}(x, s) = -\hat{Y}_{2x}(x, s),\\ 
\hat{Y}_2(0,s)=-c_2\hat{z}(0,s).
\end{array}
\right.
\end{equation}
It should be noted that although \eqref{28} is a direct copy of \eqref{27}, the initial value of \eqref{28} can be set while that of \eqref{27} is unknown. 
Define the observer error 
\begin{equation}\label{29}
\tilde{z}(x,s)=z(x,s)-\hat{z}(x,s), \;
\tilde{Y}_2(x,s)=Y_2(x,s)-\hat{Y}_2(x,s),
\end{equation}
which satisfy
\begin{equation}\label{30}
\left\{
\begin{array}{l}
\tilde{z}_{s s}(x, s) = \tilde{z}_{x x}(x, s), \\
\tilde{z}_x(0, s) = 0, \\
\tilde{z}(1, s) = -\tilde{Y}_2(1,s), \\
\tilde{Y}_{2s}(x, s) = -\tilde{Y}_{2x}(x,s),\\
\tilde{Y}_2(0,s) = -c_2\tilde{z}(0,s).
\end{array}
\right.
\end{equation}
We consider \eqref{30} in the state space 
$
\mathcal{H}_4 = \{ H^1(0,1) \times L^2(0,1) \times H^1(0,1) | f(1)=-h(1), h(0)=-c_2 f(0) \}
$.
\begin{lemma}\label{le4}
Suppose that $0<c_2<1$ and $t>\tau$.
For any initial value 
$( \tilde{z}(\cdot,0), \tilde{z}_s(\cdot,0),$ $\tilde{Y}_{2}(\cdot,0) ) \in \mathcal{H}_4$, 
there exists a unique solution $( \tilde{z}, \tilde{z}_s, \tilde{Y}_{2}) \in C([0,t-\tau];\mathcal{H}_4)$ to \eqref{30} such that 
\begin{equation}\label{31}
\| ( \tilde{z}(\cdot,s), \tilde{z}_s(\cdot,s), \tilde{Y}_{2}(\cdot,s) ) \|_{\mathcal{H}_4} 
\leq 
M_2 e^{-\mu_2 s},
\end{equation}
for some $M_2, \mu_2 >0$.
\end{lemma}
The proof of Lemma \ref{le4} is similar to that of \cite[Theorem 3.1]{2}, so we omit it here. 

In what follows, we propose an adaptive observer for \eqref{17} based on a known measurement $y_d(s)$,
which is defined as 
\begin{equation}\label{32}
y_d(s) := e(s+\tau)-\hat{z}(0,s) = g_1(0)\eta(s)+\tilde{z}(0,s).
\end{equation}
Therefore, the system \eqref{17} can be rewritten as
\begin{equation}\label{33}
\left\{
\begin{array}{l}
\dot{\eta}(s)=S_\eta \eta(s),\\
y_d(s)=g_1(0)\eta(s)+\tilde{z}(0,s).
\end{array}
\right.
\end{equation}
\begin{lemma}\label{le5}
The pair $(S_\eta,g_1(0))$ is observable for every $\omega\in(0,\infty)$.
\end{lemma}
Lemma \ref{le5} guarantees that there exists a coordinate transformation 
\begin{equation}\label{47}
d(s) = T \eta(s) = (d_1(s), d_2(s))^{\top} \in \mathbb{R}^2    
\end{equation}
such that $d(s)$ satisfies
\begin{equation}\label{34}
\left
\{
\begin{array}{l}
\dot{d}(s)=S_c(\theta) d(s), \\
y_d(s)=\gamma_\eta d(s)+\tilde{z}(0,s),
\end{array}
\right.
\end{equation}
where 
\begin{equation}\label{35}
S_c(\theta)=
\left[
\begin{array}{cc}
0       & 1  \\
-\theta & 0
\end{array}
\right], \;
\theta=\omega^2, \;
\gamma_\eta=[1,0],
\end{equation}
and $ T\in\mathbb{C}^{2 \times 2}$ is a nonsingular matrix determined by 
\begin{equation}\label{36}
S_c(\theta)=T S_\eta T^{-1}, \;
\gamma_\eta =g_1(0)T^{-1}.
\end{equation}
Inspired by \cite{3}, we design the following adaptive observer for \eqref{34} as
\begin{equation}\label{37}
\left\{
\begin{array}{l}
\dot{\xi}(s)=- \iota\xi(s)-y_d(s), \\	
\dot{\hat{\chi}}_1(s)
=\hat{\varphi}(s)+\iota y_d(s)+\hat{\theta}(s) \xi(s)+k_0(y_d(s)-\hat{\chi}_1(s)), \\
\dot{\hat{\varphi}}(s)
=-\iota \hat{\varphi}(s)-\iota^2 y_d(s), \\
\dot{\hat{\theta}}(s)=k_1 \xi(s)(y_d(s)-\hat{\chi}_1(s)),
\end{array}
\right.
\end{equation}
where $\iota>0, k_0>\frac{1}{4 \iota}, k_1>0$.
Based on \eqref{37}, we set the estimate $\hat{d}(s)$ of $d(s)$ as 
\begin{equation}\label{89}
\hat{d}(s)=(\hat{d}_1(s), \hat{d}_2(s))^\top,
\end{equation}
where
\begin{equation}\label{91}
\left\{
\begin{array}{l}
\hat{d}_1(s)=\hat{\chi}_1(s), \\
\hat{d}_2(s)=\hat{\varphi}(s)+\xi(s) \hat{\theta}(s)+\iota \hat{\chi}_1(s).
\end{array}
\right.
\end{equation}
Similar to \cite[Theorem 2.1]{4}, we can show the following result.
\begin{lemma}\label{le6}
For any initial state 
$(\xi(0), \hat{\chi}_1(0), \hat{\varphi}(0), \hat{\theta}(0))^ {\top} \in \mathbb{R}^4$, we have
\begin{equation}\label{38}
\lim _{s \rightarrow \infty}
|\hat{\theta}(s)-\theta|=0 \; {\rm{exponentially}}, \;
\lim_{s \rightarrow \infty}
\|\hat{d}(s)-d(s)\|_{\mathbb{R}^2}=0 \; {\rm{exponentially}}.	
\end{equation}
\end{lemma}
Set two new functions $f_1(x)$ and $f_2(x)$, which satisfy
\begin{equation}\label{39}
\left\{
\begin{array}{l}
f_1^{\prime\prime}(x) = f_1(x)S_c^2(\theta),\\
f_1^{\prime}(0)=0,\; f_1(0) = \gamma_\eta,\\
f_2^\prime(x) = -f_2(x)S_c(\theta),\; f_2(0) = -c_2\gamma_\eta.
\end{array}
\right.
\end{equation}
Together with \eqref{22}, \eqref{36}, and \eqref{39}, we have 
\begin{equation}\label{40}
f_1(x)=g_1(x)T^{-1}, \; 
f_2(x)=g_2(x)T^{-1}.
\end{equation}
By Lemma \ref{le3}, together with \eqref{40}, we easily obtain the following result.
\begin{lemma}\label{le7}
There exists a unique solution $f_1\in C^2( [0,1];\mathbb{R}^{1 \times 2} )$ and 
$f_2 \in C^1([0,1]; $ $\mathbb{R}^{1\times 2})$ to the initial value problem \eqref{39}.
\end{lemma}
According to \eqref{21}, \eqref{47}, \eqref{34} and \eqref{40}, we set $\hat{\varepsilon}(x,s)$ and $\hat{\varepsilon}_s(x,s)$ as 
the estimated values of $\varepsilon(x,s)$ and $\varepsilon_s(x,s)$, which satisfy
\begin{equation}\label{41}
\left\{
\begin{array}{l}
\hat{\varepsilon}(x,s) = \hat{z}(x,s) + f_1(x,\hat{\theta}(s)) \hat{d}(s), \\
\hat{\varepsilon}_s(x,s) = \hat{z}_s(x,s) + f_1(x,\hat{\theta}(s)) S_c(\hat{\theta}(s)) \hat{d}(s).
\end{array}
\right.
\end{equation}
According to \eqref{29} and \eqref{41}, we set the errors between $(\varepsilon,\varepsilon_s)$ and $(\hat{\varepsilon},\hat{\varepsilon}_s)$ as
\begin{equation}\label{90}
\left\{
\begin{array}{l}
 \tilde{\varepsilon}(x,s)=\varepsilon(x,s)-\hat{\varepsilon}(x,s), \\
 \tilde{\varepsilon}_s(x,s)=\varepsilon_s(x,s)-\hat{\varepsilon}_s(x,s),
\end{array}
\right.
\end{equation}
where 
\begin{equation}\label{42}
\left\{
\begin{array}{l}
\tilde{\varepsilon}(x,s) = \tilde{z}(x,s) + f_1(x,\theta) d(s) - f_1(x,\hat{\theta}(s)) \hat{d}(s),  \\
\tilde{\varepsilon}_s(x,s) = 
\tilde{z}_s(x,s) + f_1(x,\theta) S_c(\theta) d(s) - f_1(x,\hat{\theta}) S_c(\hat{\theta}(s)) \hat{d}(s). 
\end{array}
\right.
\end{equation}
We have the following result, the proof of which will be given in subsection 3.1. 
\begin{theorem}\label{th1}
Suppose that $t>\tau$.
For any initial condition 
$( \tilde{\varepsilon} (\cdot,0), \tilde{\varepsilon}_s (\cdot,0) )\in \mathcal{H}_3$, 
there exists a unique solution $(\tilde{\varepsilon}, \tilde{\varepsilon}_s)\in C([0,t-\tau];\mathcal{H}_3)$ 
to the error system $(\tilde{\varepsilon}(x,s),\tilde{\varepsilon}_s(x,s))$ such that 
\begin{equation}\label{43}
\|( \tilde{\varepsilon}(\cdot,s), \tilde{\varepsilon}_s (\cdot,s) )\|_{\mathcal{H}_3} 
\leq M_3 e^{-\mu_3 s},
\end{equation}
for some $M_3,\mu_3>0$. 
\end{theorem}

\subsubsection{Predictor design}
In this part, we design predictors to estimate the states of the coupled system \eqref{9} and \eqref{34} for $s\in [t-\tau,t], \; t>\tau$.
First, we propose the following predictor for \eqref{34}:
\begin{equation}\label{68}
\left\{
\begin{array}{l}
\dot{\hat{d}} ^t (s) = S_c(\hat{\theta}(t-\tau)) \hat{d}^t (s),
\; t-\tau < s < t, t>\tau,\\
\hat{d} ^t (t-\tau)=\hat{d}(t-\tau).   
\end{array}
\right.
\end{equation}
According to \eqref{22}, \eqref{47}, and \eqref{40}, we have
\begin{equation}\label{48}
\gamma_\eta \eta(s) = -( f_1(1,\theta) + f_2(1,\theta)) d(s).
\end{equation}
Based on \eqref{37}, \eqref{68} and \eqref{48}, we define the estimated value $D^t (s)$ of $\gamma_\eta \eta(s)$ as 
\begin{equation}\label{67}
D^t (s) = -( f_1(1,\hat{\theta}(t-\tau) + f_2(1,\hat{\theta}(t-\tau))) \hat{d}^t (s).
\end{equation}
Combining with \eqref{67}, we design the following predictor for \eqref{9}:
\begin{equation}\label{49}
\left
\{
\begin{array}{l}
\hat{\varepsilon}^t_{s s}(x, s)
=\hat{\varepsilon}^t_{x x}(x, s), 
\; 0<x<1, t-\tau<s<t, t>\tau, \\
\hat{\varepsilon}^t_x(0, s)
=-q \hat{\varepsilon}^t(0, s), 
\; t-\tau \leq s \leq t, t>\tau, \\
\hat{\varepsilon}^t(1, s)
=u(s) - D^t (s), 
\; t-\tau \leq s \leq t, t>\tau, \\
\hat{\varepsilon}^t(x, t-\tau) = \hat{\varepsilon}(x, t-\tau), \;
\hat{\varepsilon}^t_s(x, t-\tau) = \hat{\varepsilon}_s(x, t-\tau), 
\; 0 \leq x \leq 1, t>\tau .
\end{array}
\right.
\end{equation}
\begin{remark}\label{re2}
The design idea of the predictor \eqref{68} and \eqref{49} is stated as:
First, we replace $\theta$ with $\hat{\theta}(t-\tau)$ to estimate the disturbance $d(t)$ for \eqref{68}.
Then, we use the final states of observers \eqref{89} and \eqref{41} as initial states of \eqref{68} and \eqref{49}.
In addition, predictors \eqref{68} and \eqref{49} replicates the dynamics of the systems \eqref{34} and \eqref{9}.
\end{remark}
Setting the predictor error $\tilde{\varepsilon}^t(x,s)$ as
\begin{equation}\label{50}
\tilde{\varepsilon}^t(x, s)
=\varepsilon(x, s)-\hat{\varepsilon}^t(x, s), 
\end{equation}
together with \eqref{9} and \eqref{49}, we have                                                           
\begin{equation}\label{51}
\left
\{
\begin{array}{l}
\tilde{\varepsilon}^t_{s s}(x, s)
=\tilde{\varepsilon}^t_{x x}(x, s), 
\; 0<x<1, t-\tau<s<t, t>\tau, \\
\tilde{\varepsilon}^t_x(0, s)
=-q \tilde{\varepsilon}^t(0, s), 
\\
\tilde{\varepsilon}^t(1, s)
= \Delta^t(s),
\\
\tilde{\varepsilon}^t(x, t-\tau) = \tilde{\varepsilon}(x, t-\tau), \;
\tilde{\varepsilon}^t_s(x, t-\tau) = \tilde{\varepsilon}_s(x, t-\tau), 
\end{array}
\right.
\end{equation}
where
\begin{equation}\label{52}
\Delta^t(s)=
(f_1(1,\theta)+f_2(1,\theta)) d(s)
+ D^t (s).
\end{equation}

We have the following result.
\begin{theorem}\label{th2}
Suppose that $t>\tau$, for any initial value 
$ ( \tilde{\varepsilon}(\cdot, t-\tau), \tilde{\varepsilon}_s(\cdot, t-\tau)) \in \mathcal{H}_1 $,
there exists a unique solution 
$( \tilde{\varepsilon}^t, \tilde{\varepsilon}^t_s) \in C([t-\tau,t] ; \mathcal{H}_1)$ 
to \eqref{51} such that
\begin{equation}\label{55}
 \| ( \tilde{\varepsilon}^t(\cdot,t), \tilde{\varepsilon}^t_s(\cdot,t) )\| _{\mathcal{H}_1}
\leq
M_{5} e^{-\mu_{5} t},
\end{equation}
for some $M_5,\mu_5>0$.
\end{theorem}

\subsection{Error-based feedback controller}
In this section, we propose the output feedback controller for \eqref{1} only based on $e(t)$ and prove the exponential stability of the final closed-loop system. 

According to \eqref{67} and \eqref{49}, we substitute the unknown variables of \eqref{19} by the corresponding estimates and obtain the final feedback controller
\begin{equation}\label{69}
u(t)=
\left\{
\begin{array}{l}
0, 
\; 0 \leq t \leq \tau, \\
-\dfrac{1}{c_0\!+\!q}\!
[\hat{\varepsilon}^t_x(1, t)\!+\!c_1 \!\hat{\varepsilon}^t_s(1, t)]
\!-\!\!\int_0^1 e^{q(1\!-\!h)}
[q \hat{\varepsilon}^t(h, t) \!+\! c_1 \!\hat{\varepsilon}^t_s(h, t)] d h\!+\! D^t(s), t\!>\!\tau.
\end{array}
\right.
\end{equation}
Together with \eqref{1}, \eqref{2}, \eqref{3}, \eqref{20}, \eqref{28}, \eqref{32}, \eqref{37}, \eqref{41}, \eqref{68}, \eqref{67}, \eqref{49}, and \eqref{69}, the final closed-loop system is described as 
\begin{equation}\label{70}
\left
\{
\begin{array}{l}
w_{t t}(x, t)=w_{x x}(x, t)+p_1(x)v(t), \; t>0,\\
w_x(0, t)=-q w(0, t)+p_2v(t), \;
w(1, t)=u(t)+p_3v(t), \\
\dot{v}(t)=Sv(t),\; t>0, \\
e(t)=y_p(t)-y_{\rm{ref}}(t)=w(0, t-\tau)-p_4v(t), \; t>0, \\
Y_{1s}(x,s)=-Y_{1x}(x,s), \; 0<s<t-\tau, \\
Y_1(0,s)=-c_2e(s+\tau),\; 0<c_2<1, \\ 
\hat{z}_{s s}(x, s) = \hat{z}_{x x}(x, s), 
\; 0<s<t-\tau, t>\tau, \\
\hat{z}_x(0, s)=-qe(s+\tau), \;
\hat{z}(1, s) = u(s)-(\hat{Y}_2(1,s)-Y_1(1,s)), \\
\hat{Y}_{2s}(x, s) = -\hat{Y}_{2x}(x, s), 
\; 0<s<t-\tau, t>\tau,\\ 
\hat{Y}_2(0,s)=-c_2\hat{z}(0,s), \\
y_d(s) = e(s+\tau)-\hat{z}(0,s), 
\; 0<s<t-\tau, t>\tau, \\
\dot{\xi}(s)=- \iota\xi(s)-y_d(s), \\	
\dot{\hat{\chi}}_1(s)
=\hat{\varphi}(s)+\iota y_d(s)+\hat{\theta}(s) \xi(s)+k_0(y_d(s)-\hat{\chi}_1(s)), \\
\dot{\hat{\varphi}}(s)
=-\iota \hat{\varphi}(s)-\iota^2 y_d(s), \\
\dot{\hat{\theta}}(s)=k_1 \xi(s)(y_d(s)-\hat{\chi}_1(s)),\\
\hat{d}_1(s)=\hat{\chi}_1(s), 
\;
\hat{d}_2(s)=\hat{\varphi}(s)+\xi(s) \hat{\theta}(s)+\iota \hat{\chi}_1(s), \\
\hat{\varepsilon}_{s s}^t(x, s) = \hat{\varepsilon}^t_{x x}(x, s), 
\; t-\tau<s<t, t>\tau, \\
\hat{\varepsilon}^t_x(0, s) = -q \hat{\varepsilon}^t(0, s), 
\;
\hat{\varepsilon}^t(1, s) = u(s) - D^t(s), \\
\hat{\varepsilon}^t(x, t-\tau) = \hat{z}(x, t-\tau) + f_1(x,\hat{\theta}(t-\tau)) \hat{d}(t-\tau), \\
\hat{\varepsilon}^t_s(x, t-\tau) = \hat{z}_s(x, t-\tau) 
+ f_1(x,\hat{\theta}(t-\tau)) S_c(\hat{\theta}(t-\tau)) \hat{d}(t-\tau), \\
\dot{\hat{d}} ^t(s) = S_c(\hat{\theta}(t-\tau)) \hat{d}^t(s),
\; t-\tau < s < t, t>\tau,\\
\hat{d}^t(t-\tau)=\hat{d}(t-\tau),\\ 
D^t(s) = -( f_1(1,\hat{\theta}(t-\tau) + f_2(1,\hat{\theta}(t-\tau))) \hat{d}^t(s),
\; t-\tau < s < t, t>\tau.
\end{array}
\right.
\end{equation}
We consider the closed-loop system \eqref{70} in the state space 
$
\mathcal{H}_5=\mathcal{H}_2 \times H^{1}(0,1) \times \mathcal{H}_4 \times \mathbb{R}^4 \times \mathcal{H}_2.
$ We have the following result.
\begin{theorem}\label{th3}
Suppose that the parameters $q, c_0, c_1>0$, $0<c_2<1$, $\iota,k_1>0$, $k_0>\frac{1}{4\iota}$ and $t>\tau$.
For any initial value 
$( 
w_0(\cdot), w_1(\cdot), v_0, Y_1(\cdot,0), \hat{z}(\cdot,0), 
$
$
\hat{z}_s(\cdot,0), \hat{Y}_2(\cdot,0),$ $
\xi(\cdot), \hat{\chi}_1(\cdot), \hat{\varphi}(\cdot), \hat{\theta}(\cdot) 
)$ 
$
\in
\mathcal{H}_2 \times H^1(0,1) \times \mathcal{H}_4 \times \mathbb{R}^4
$,
there exists a unique solution to \eqref{70} such that
$(w, w_t, v) \in C([0,\infty) ; \mathcal{H}_2)$, 
$
( Y_1, \hat{z}, \hat{z}_s, \hat{Y}_2, \xi, \hat{\chi}_1,
$
$\hat{\varphi}, \hat{\theta}) 
\in C([0,t-\tau];H^1(0,1) \times \mathcal{H}_4 \times \mathbb{R}^4 )$,
and 
$( \hat{\varepsilon}^t, \hat{\varepsilon}^t_s, \hat{d}^t) 
\in C([t-\tau,t] \times [\tau,\infty); \mathcal{H}_2)$.
Moreover, the tracking error $e(t)$ is exponentially convergent to $0$, that is,
\begin{equation}\label{71}
\lim_{t\to\infty}|e(t)| = 0,\; \rm{exponentially}.
\end{equation}
\end{theorem}

\section{Proof of main results}

\subsection{Proof of Theorem \ref{th1}}
According to Lemma \ref{le6}, we get that there exist positive constants $M_4$ and $\mu_4$ such that $|\hat{\theta}(s)-\theta|< M_4 e^{-\mu_4 s}$.
Thus, we obtain
\begin{equation}\label{44}
\theta-M_4 \leq |\hat{\theta}(s)| \leq \theta+M_4.
\end{equation}

From \eqref{35} and \eqref{39}, we know that $f_1(x,\theta),f_2(x,\theta),S_c(\theta)$ are all continuously differentiable with respect to $\theta$, so they are all Lipschitz continuous functions over $[ \theta-M_4,\theta+M_4]$.
Therefore, we demonstrate that there exists $L>0$, such that  
\begin{equation}\label{45}
\left\{
\begin{array}{ll}
\|f_1(x,\hat{\theta})-f_1(x,\theta)\|
\leq L|\hat{\theta}(s)-\theta|
\leq L M_4 e^{-\mu_4 s},
& \forall x\in[0,1],\\\|f_1^{\prime}(1,\hat{\theta})-f_1^{\prime}(1,\theta)\|
\leq L|\hat{\theta}(s)-\theta|
\leq L M_4 e^{-\mu_4 s}, \\
\|f_2(1,\hat{\theta})-f_2(1,\theta)\|
\leq L|\hat{\theta}(s)-\theta|
\leq L M_4 e^{-\mu_4 s},\\
\|S_c(\theta)-S_c(\hat{\theta})\|
\leq L|\hat{\theta}(s)-\theta|
\leq L M_4 e^{-\mu_4 s},\\
\|e^{S_c(\theta)s}-e^{S_c(\hat{\theta})s}\|
\leq L|\hat{\theta}(s)-\theta|
\leq L M_4 e^{-\mu_4 s},
\end{array}
\right.
\end{equation}
for $s\in[0,t-\tau]$.
Hence, combining Lemma \ref{le6}, \eqref{45} with \cite[Lemma 2.5]{3}, we have for $x\in[0,1]$
\begin{equation}\label{46}
\left\{
\begin{array}{l}
\lim_{s \to \infty}
|f_1(x,\hat{\theta})\hat{d}(s)-f_1(x,\theta)d(s)|
= 0,
\; \rm{exponentially}, \\
\lim_{s \to \infty}
|f_2(1,\hat{\theta})\hat{d}(s)-f_2(1,\theta)d(s)|
= 0,
\; \rm{exponentially},\\
\lim_{s \to \infty}
|f_1(x,\hat{\theta})S_c(\hat{\theta})\hat{d}(s)-f_1(x,\theta)S_c(\theta)d(s)|
= 0,
\; \rm{exponentially}.
\end{array}
\right.
\end{equation}
According to \eqref{31}, \eqref{42}, and \eqref{46}, we finally obtain \eqref{43}.

\subsection{Proof of Theorem \ref{th2}}
According to \eqref{52}, we have 
{\small{
\begin{align}\label{53}
& \Delta^{t \prime}(s) \cr
= & (f_1(1,\theta)+f_2(1,\theta)) S_c(\theta)  d(s) \cr
&-( f_1 (1,\hat{\theta}(t-\tau)) + f_2 (1,\hat{\theta}(t-\tau)) ) S_c(\hat{\theta}(t-\tau)) \hat{d}^t(s) \cr
= & (f_1(1,\theta)+f_2(1,\theta)) S_c(\theta) e^{ S_c(\theta)(s-t+\tau) } d(t-\tau) \cr
& - \!(f_1 (1,\hat{\theta}(t-\tau)) \! + \! f_2 (1,\hat{\theta}(t-\tau)) ) S_c(\hat{\theta}(t-\tau)) 
e^{ S_c(\hat{\theta}(t-\tau))(s-t+\tau) }
\hat{d}(t-\tau),
\end{align}
}}
where $\Delta^{t \prime}(s):=\dfrac{{\rm d}\Delta^{t}(s)}{{\rm d} s}$.
According to \cite[Lemma 2.5]{3}, Lemma \ref{le6}, \eqref{45}, \eqref{46}, and \eqref{53}, 
we obtain 
\begin{equation}\label{54}
\begin{array}{ll}
\lim_{ t\to\infty } \Delta^t (s) = 0, &  {\rm{ exponentially }}, \\
\lim_{ t\to\infty } \Delta^{t \prime} (s) = 0,  & {\rm{ exponentially }},
\end{array}
\end{equation}
for $s\in [t-\tau,t]$.
We divide the rest of the proof into two main parts:
the well-posedness and the exponential stability of the error system \eqref{51}. 

\noindent{\bf{Step 1.} We prove the well-posedness of \eqref{51}, $\forall s\in[t-\tau,t]$.}

We consider \eqref{51} in the state space $\mathcal{H}_1$.
Inspired by the proof of \cite[Lemma 6.1]{2}, we propose the following transformation to fix this problem:
\begin{equation}\label{56}
\epsilon^t(x, s) = \tilde{\varepsilon}^t(x, s)- \Delta^t_e(s-1+x),
\end{equation}
where 
\begin{equation}\label{57}
\Delta^t_e(h)=
\left \{
\begin{array}{ll}
\Delta^t(t-\tau), & 
h \leq t-\tau, \\
\Delta^t(h), & 
h \geq t-\tau.
\end{array}
\right.
\end{equation}
Hence, by \eqref{56} and \eqref{57}, the system \eqref{51} is transformed into 
\begin{equation}\label{58}
\left
\{
\begin{array}{l}
\epsilon^t_{s s} (x, s)
=\epsilon^t_{x x} (x, s),
\; 0<x<1, t-\tau<s<t, t>\tau, \\
\epsilon^t_x (0, s)
=-q \epsilon^t (0, s) - q\Delta^t_e(s-1) - \Delta ^{t \prime}_{e}(s-1),
\\
\epsilon^t(1, s)
=0,
\\
\epsilon^t (x, t-\tau)
=\tilde{\varepsilon}(x, t-\tau) - \Delta^t(t-\tau), \;
\epsilon^t_s(x, t-\tau) = \tilde{\varepsilon}_s(x, t-\tau). 
\end{array}
\right.
\end{equation}
Now, our problem reduces to proving the well-posedness of \eqref{58}.
This system is considered in the space $\mathcal{H}_1$ with the inner product 
{\small{
\begin{align}\label{59}
&\langle (\phi_1, \psi_1 ), (\phi_2, \psi_2 ) \rangle_{\mathcal{H}_1} \cr
=&\int_0^1 
[\phi_1^{\prime}(x) \overline{ \phi_2^{\prime}(x) } + \psi_1(x) \overline{ \psi_2(x) } ] d x
+\phi_1(0) \overline{\phi_2(0)} ,
\; \forall \phi_j,\psi_j\in \mathcal{H}_1, j=1,2.
\end{align}
}}
By deforming \eqref{58}$_2$, the abstract evolutionary equation of the system \eqref{58} is defined as
\begin{equation}\label{60}
\dfrac{ d }{ d t }
\left[
\begin{array}{c}
\epsilon^t (\cdot, s)\\
\epsilon^t_s (\cdot, s)
\end{array}
\right]
=\mathbb{A}
\left[
\begin{array}{c}
\epsilon^t (\cdot, s)\\
\epsilon^t_s (\cdot, s)
\end{array}
\right]
+\mathbb{B}[(q+1)\epsilon^t (0,s) ] 
+\mathbb{B} [ q\Delta^t_e(s-1) + \Delta ^{t \prime}_{e}(s-1) ] ,
\end{equation}
where the operators 
$
\mathbb{A}: 
\mathcal{D}(\mathbb{A}) \subset \mathcal{H}_1 
\rightarrow \mathcal{H}_1, 
\mathbb{B}: 
\mathbb{C} 
\rightarrow \mathcal{D} (\mathbb{A}^* )^{\prime}
$
are given by
\begin{equation}\label{61}
\left\{
\begin{array}{l}
\mathbb{A}(\phi, \psi)
=(\psi, \phi^{\prime \prime}), \; \forall(\phi, \psi) 
\in \mathcal{D}(\mathbb{A}), \\
\mathcal{D}(\mathbb{A})
=\{
(\phi, \psi) 
\in H^2(0,1) \times H^1_R(0,1) 
\mid \phi^{\prime}(0)=\phi(0)
\}, \\
\mathbb{B}
=(0,\delta(x) ).
\end{array}
\right.
\end{equation}
According to \cite{6}, we have 
$\mathcal{D}(\mathbb{A}^*) = \mathcal{D}(\mathbb{A}), \; \mathbb{A}^*=-\mathbb{A}$
and
$\mathbb{B}^*(\phi, \psi)=\psi(0), \; \forall(\phi, \psi) \in \mathcal{D} (\mathbb{A}^*)$.
Hence, we demonstrate that $\mathbb{B}$ is admissible for $e^{\mathbb{A}t}$.
Together with \cite[Proposition 1.1]{7}, we finish the proof of the well-posedness of \eqref{58}.
Combining this result with \eqref{56} and \eqref{57}, we can further prove that the system \eqref{51} is well-posed.

\noindent{\bf{Step 2.} We prove the exponential stability of \eqref{51}.}

In this part, we still utilize \eqref{58} to complete this proof.
In addition to \eqref{60} and \eqref{61}, there are another abstract evolutionary equation for \eqref{58}:
\begin{equation}\label{62}
\dfrac{ d }{ d t }
\left[
\begin{array}{c}
\epsilon^t (\cdot, s)\\
\epsilon^t_s (\cdot, s)
\end{array}
\right]
=\mathcal{A}
\left[
\begin{array}{c}
\epsilon^t (\cdot, s) \\
\epsilon^t_s (\cdot, s)
\end{array}
\right]
+\mathcal{B}[ q\Delta^t_e(s-1) + \Delta ^{t \prime}_{e}(s-1) ],
\end{equation}
where 
$\mathcal{A}: \mathcal{D}(\mathcal{A}) \subset \mathcal{H}_1 
\rightarrow \mathcal{H}_1, \;
\mathcal{B}: \mathbb{C} 
\rightarrow \mathcal{D} (\mathcal{A}^* )^{\prime}$
are defined as 
\begin{equation}\label{63}
\left\{
\begin{array}{l}
\mathcal{A}(\phi, \psi)
=(\psi, \phi^{\prime \prime}), \; \forall(\phi, \psi) 
\in \mathcal{D}(\mathcal{A}), \\
\mathcal{D}(\mathcal{A})
=\{
(\phi, \psi) 
\in H^2(0,1) \times H_R^1(0,1) 
\mid \phi^{\prime}(0)= -q\phi(0)
\}, \\
\mathcal{B}
=(0,\delta(x)).
\end{array}
\right.
\end{equation}
According to \cite[Lemma 2.1]{8}, we demonstrate that $\mathcal{A}$ generates a $C_0$ semigroup,
further, we prove that $\mathcal{B}$ is admissible for $e^{\mathcal{A} t}$ by the well-posedness of \eqref{58}.
In what follows, we focus on the moment $s=t$, and the system \eqref{58} can be written as
\begin{align} \label{64}
& ( \epsilon^t(x, t), \epsilon^t_s(x, t) )\cr
= & ( \epsilon (x, t-\tau),\epsilon_s(x, t-\tau) ) e^{ \mathcal{A} \tau} 
+ \int_{ t-\tau }^{ t } 
e^{ \mathcal{A} (t-h)} 
\mathcal{B} [ q \Delta^t_e(h-1) + \Delta ^{t \prime}_{e}(h-1) ] dh \cr
= &( \epsilon (x, t-\tau),\epsilon_s(x, t-\tau) ) e^{ \mathcal{A} \tau} \cr
& + \int_{ 0 }^{ \tau } 
e^{ \mathcal{A} (\tau-\sigma)} 
\mathcal{B} 
[q \Delta^t_e(\sigma+t-\tau-1) + \Delta ^{t \prime}_{e}(\sigma+t-\tau-1)] d\sigma,
\end{align}
where $\sigma=h-(t-\tau), \; \sigma\in[0,\tau]$.
According to Theorem \ref{th1}, \eqref{54}, and \cite[Proposition 2.5]{5}, we can scale \eqref{64} as
\begin{align}\label{65}
& \| (\epsilon^t (\cdot, t), \epsilon^t_s(\cdot, t)) \| _{\mathcal{H}_1} \cr
\leq & M_6 e^{\mu_6 \tau}
\| ((\epsilon^t (\cdot, t-\tau), \epsilon^t_s(\cdot, t-\tau) ) \| _{\mathcal{H}_1} \cr
& + \mathcal{L}_{\mathcal{B}} |q \Delta^t_e(\sigma+t-\tau-1) + \Delta ^{t \prime}_{e}(\sigma+t-\tau-1)| \cr
\leq & M_7 e^{-\mu_7 t}, 
\end{align}
for some $ M_6,M_7,\mu_6,\mu_7,\mathcal{L}_{\mathcal{B}} > 0$.
Combining \eqref{54}, \eqref{56} with \eqref{65}, we get
\begin{align}\label{66}
& \dfrac{1}{2} \int_0^1 
\left[|\tilde{\varepsilon}^t_x(x, t)|^2 
+|\tilde{\varepsilon}^t_s(x, t)|^2 \right] d x \cr
\leq &
\| (\epsilon^t(\cdot,t), \epsilon^t_s(\cdot, t) )\| _{\mathcal{H}_1}^2
+2 \int_0^1 |\Delta ^{t \prime}_{e} (t-1+x)|^2 d x \cr
= & \| (\epsilon^t(\cdot,t), \epsilon^t_s(\cdot, t) )\| _{\mathcal{H}_1}^2
+2 \int_{t-1}^t| \Delta ^{t \prime}_{e} (h) |^2 d h \cr
\leq & \frac{1}{2} M_5^2 e^{-2\mu_5 t}.
\end{align}
The proof is complete.
\subsection{Proof of Theorem \ref{th3}}
First, we divide the time interval into two parts: $t\in[0,\tau]$ and $t\in[\tau,\infty]$.

\noindent{\bf{Step 1.} We prove the well-posedness of \eqref{9} under the controller \eqref{69} for $0<t<\tau$.}

According to \eqref{48} and \eqref{69}, when $t\in[0,\tau]$, the system \eqref{9} is described as 
\begin{equation}\label{72}
\left
\{
\begin{array}{l}
\varepsilon_{tt}(x,t)
=\varepsilon_{x x}(x, t), \\
\varepsilon_x(0, t)
=-q\varepsilon(0,t),\\
\varepsilon(1, t)
= -(f_1(1,\theta)+f_2(1,\theta)) d(t).
\end{array}
\right.
\end{equation}
Together with \eqref{34}$_1$, we know that $d(t)$ is continuously differentiable.
Similar to the well-posedness part of the proof for Theorem \ref{th2}, we demonstrate that the system \eqref{9} is well-posed on $[0,\tau]$ and obtain
\begin{equation}\label{74}
\| ( \varepsilon(\cdot, \tau),\varepsilon_t(\cdot,\tau) ) \|_{\mathcal{H}_1}
\leq 
M_6 e^{\mu_6 \tau} \| ( \varepsilon(\cdot, 0),\varepsilon_t(\cdot,0 )) \|_{\mathcal{H}_1} 
+M_{7},
\end{equation}
for some $M_6, \mu_6, M_7>0$.

\noindent{\bf{Step 2. We prove the well-posedness and stability of \eqref{9} under the controller \eqref{69} for $t>\tau$.}}

When $s=t$, together with the controller \eqref{69}, the system \eqref{9} is described as
{\small{
\begin{equation}\label{75}
\left
\{
\begin{array}{l}
\varepsilon_{tt}(x,t)
=\varepsilon_{x x}(x, t), \\
\varepsilon_x(0, t)
=-q\varepsilon(0,t),\\
\varepsilon(1, t)
=-\dfrac{1}{c_0+q}
[\varepsilon_x(1, t)-\tilde{\varepsilon}^t_x(1, t)
+c_1 \varepsilon_t(1, t)-c_1\tilde{\varepsilon}^t_s(1, t)] \\
\hspace{4em} -\int_0^1 e^{q(1-h)}
[q ( \varepsilon(h,t)-\tilde{\varepsilon}^t(h, t) )
+ c_1 (\varepsilon_t(h, t)-\tilde{\varepsilon}^t_s(h, t) )] d h
+\Delta^t(t).
\end{array}
\right.
\end{equation}
}}
According to the Volterra transformation \eqref{10}, the system \eqref{75} can be transformed into 
\begin{equation}\label{76}
\left
\{
\begin{array}{l}
\bar{\varepsilon}_{tt}(x,t)
=\bar{\varepsilon}_{x x}(x, t), \\
\bar{\varepsilon}_x(0, t)
=c_0 \bar{\varepsilon}(0,t),\\
\bar{\varepsilon}_x(1, t)
=-c_1 \bar{\varepsilon}_t (1,t)
+\tilde{\varepsilon}^t_x(1, t)+c_1 \tilde{\varepsilon}^t_s(1, t) \\
\hspace{4em} +(c_0+q) \int_0^1 e^{q(1-h)}
(q \tilde{\varepsilon}^t(h, t)
+c_1 \tilde{\varepsilon}^t_s(h, t)) d h
+(c_0+q)\Delta^t(t).
\end{array}
\right.
\end{equation}
We consider \eqref{76} in the state space $\mathcal{H}_3$ with the inner product:
{\small{
\begin{align}\label{78}
& \langle ( \phi_1, \psi_1 ),( \phi_2, \psi_2 ) \rangle_{\mathcal{H}_3} \cr
=& \int_0^1 \phi_1^{\prime}(x) \overline{\phi_2^{\prime}(x)} d x
\! + \! \psi_1(x) \overline{\psi_2(x)} d x 
\!+ \! c_0 \phi_1(0) \overline{\phi_2(0)}, 
 \forall (\phi_j, \psi_j) \in \mathcal{H}_3, j=1,2 .
\end{align}
}}
Hence, the system \eqref{76} can be written as
\begin{align}\label{79}
\frac{d}{d t}
\binom{\bar{\varepsilon}(\cdot, t)}{\bar{\varepsilon}_t(\cdot, t)}
= & \mathscr{A}
\binom{\bar{\varepsilon}(\cdot, t)}{\bar{\varepsilon}_t(\cdot, t)} \cr
& +\mathscr{B}
\left( 
(c_0+q) \int_0^1 e^{q(1-h)} (q \tilde{\varepsilon}^t(h, t)
+c_1 \tilde{\varepsilon}^t_s(h, t)) d h 
\right. \cr
& 
+\tilde{\varepsilon}^t_x(1, t)+c_1 \tilde{\varepsilon}^t_s(1, t)+(c_0+q)\Delta^t(t)
\bigg),
\end{align}
where
{\small{
\begin{equation}\label{80}
\left\{
\begin{array}{l}
\mathscr{A}(\varphi, \psi)^{\top} 
= 
(\psi, \varphi^{\prime \prime})^{\top}, \;  \forall(\varphi, \psi)^{\top} \in \mathcal{D}(\mathscr{A}), \\
\mathcal{D} ( \mathscr{A} )
=
\{
(\varphi, \psi)^{\top} \in H^1(0,1) \times L^2(0,1) 
\mid 
\varphi^{\prime}(0)=c_0 \varphi(0), 
\varphi^{\prime}(1)=-c_1 \psi(1) 
\}, \\
\mathscr{B}=
(0 \quad \delta(x-1) )^\top.
\end{array}
\right.
\end{equation}
}}
According to the proof of \cite[Theorem 5.1]{8}, we know that $\mathscr{A}$ is a generator of an exponentially stable semigroup and $\mathscr{B}$ is admissible for $e^{\mathscr{A}t}$.

Now, our problem reduces to proving that 
$e^{\frac{\alpha}{2} \cdot} \tilde{\varepsilon}^t_s(1,\cdot),
e^{\frac{\alpha}{2}\cdot} \tilde{\varepsilon}^t_x(1,\cdot)
\! \in \! L^2(\tau,\infty)$, 
where $\dfrac{\alpha}{2} \in (0,\mu_5)$.
When $s=t$, for \eqref{51}, we define
\begin{equation}\label{81}
\rho(t)=2e^{\alpha t}\int_0^1 x \tilde{\varepsilon}^t_x(x,t) \tilde{\varepsilon}^t_s(x,t) dx.
\end{equation}
According to Theorem \ref{th2}, we have 
\begin{equation}\label{82}
0 
\leq |\rho(t)| 
\leq e^{\alpha t} \int_0^1 [\tilde{\varepsilon}_x^{t2}(x, t)
+\tilde{\varepsilon}_s^{t2}(x, t)] d x
\leq M_5^2 e^{-(2\mu_5-\alpha)t},
\end{equation}
and differentiate $\rho(t)$ along with the solution of \eqref{51} gives
\begin{equation}\label{83}
\dot{\rho}(t)
=\alpha \rho(t)
+ e^{\alpha t}  [ \tilde{\varepsilon}^{t2}_x(1, t) + \tilde{\varepsilon}_s^{t2}(1, t)  ]
-e^{\alpha t} \int_0^1 
\left[ \tilde{\varepsilon}_x^{t2}(x, t)
+\tilde{\varepsilon}_s^{t2}(x,t) \right] d x.
\end{equation}
It follows from \eqref{82}, \eqref{83}, and Theorem \ref{th2} that
\begin{align}\label{84}
&\int_\tau^{\infty} 
e^{\alpha t}
\left[
\tilde{\varepsilon}_s^{t2}(1,t)
+\tilde{\varepsilon}_x^{t2}(1,t)
\right] d t \cr
\leq 
& |\rho(\infty)| + |\rho(\tau)|
+(\alpha+1) 
\int_\tau^{\infty} e^{\alpha t}
\int_0^1 [\tilde{\varepsilon}_x^{t2}(x,t)+\tilde{\varepsilon}_s^{t2}(x,t)] d x d t<+\infty .
\end{align}
According to Theorem \ref{th2}, \eqref{54}, and \cite[Lemma 2.1]{7}, we obtain
\begin{align}\label{85}
\| ( \bar{\varepsilon}(\cdot, t),\bar{\varepsilon}_s(\cdot,t) ) \|_{\mathcal{H}_3} 
\leq & 
M_8 e^{-\mu_8 (t-\tau)}
\| ( \bar{\varepsilon}(\cdot, \tau),\bar{\varepsilon}_s(\cdot,\tau) ) \|_{\mathcal{H}_3} \cr
&+\mathcal{L}_\mathscr{B}
(
\|\tilde{\varepsilon}^t_x(1, \cdot)\|_{L^2(\tau,\infty)} 
 + c_1 \|\tilde{\varepsilon}^t_s(1, \cdot)\|_{L^2(\tau,\infty)} \cr
& +  M_{9}e^{-\mu_{9}t}
+ (c_0+q) |\Delta^t(t)|
),
\end{align}
for some $M_8,\mu_8, M_9,\mu_9, \mathcal{L}_{ \mathscr{B} }>0 $.
From the inverse transformation \eqref{11}, we have
\begin{equation}\label{86}
\binom{\varepsilon(\cdot, t)}{\varepsilon_t(\cdot, t)}
=\left(
\begin{array}{cc}
(I+\mathbb{P})^{-1} & 0 \\
0 & (I+\mathbb{P})^{-1}
\end{array}
\right)
\binom{\bar{\varepsilon}(\cdot, t)}{\bar{\varepsilon}_t(\cdot, t)},
\end{equation}
where $\mathbb{P}$ is defined in \eqref{11}.
Therefore, there exists a unique solution to \eqref{75} such that $(\varepsilon(\cdot, t), \varepsilon_t(\cdot, t)) \in C(\tau, \infty ; \mathcal{H}_3)$. 
Moreover, this solution satisfies
\begin{equation}\label{87}
\|(\varepsilon(\cdot, t), \varepsilon_t(\cdot, t))^{\top}\|_{\mathcal{H}_3} 
\leq 
C_{\mathbb{P}}\|(\bar{\varepsilon}(\cdot, t), \bar{\varepsilon}_t(\cdot, t))^{\top}\|_{\mathcal{H}_3} 
\leq  M_{10} e^{-\mu_{10} t},
\end{equation}
for some $C_{\mathbb{P}},M_{10},\mu_{10}>0$.
Thus, according to the inner product \eqref{78}, we obtain that the tracking error $e(t)=\varepsilon(0,t-\tau)\to 0,\;\mbox{exponentially}$ as $t\to \infty$. 

\section{Numerical simulations}
\begin{figure}[htbp]
\centering
\subfigure[ Tracking performance $e(t)$ ]
{
\includegraphics[width=6cm]{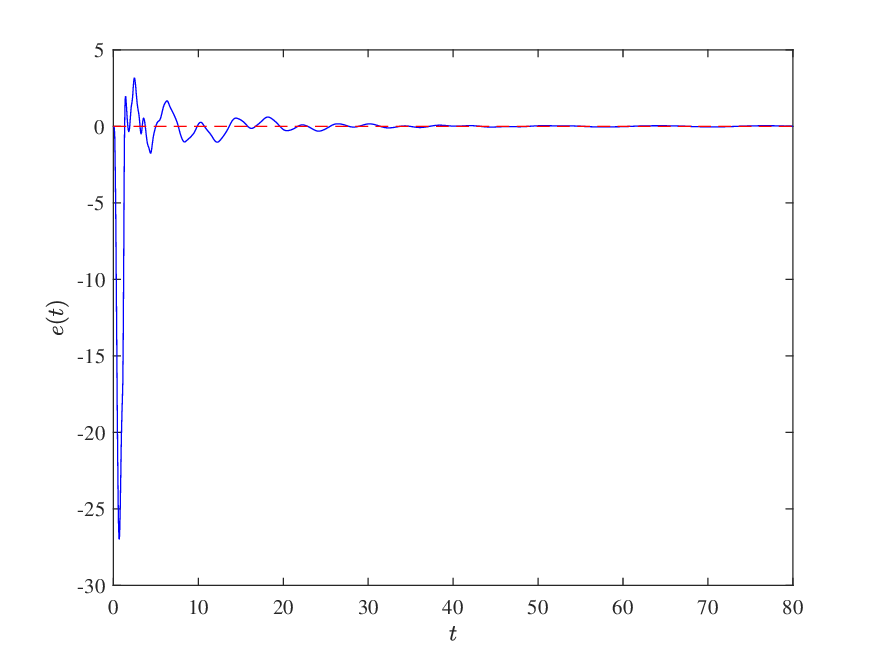}
\label{fig1}
} 
\quad 
\subfigure[ Estimate of $\theta$ ]
{
\includegraphics[width=6cm]{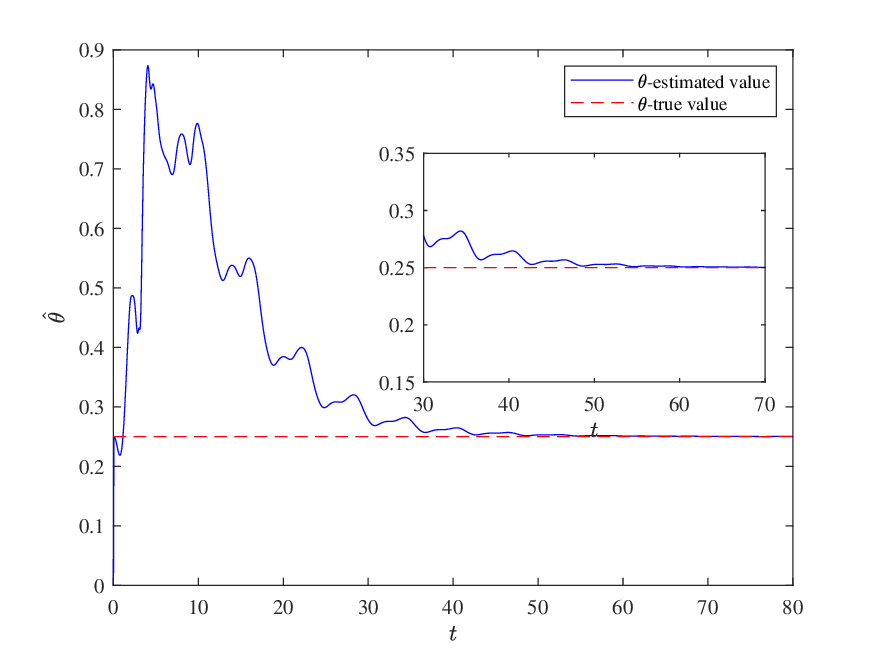}
\label{fig2}
}
\quad
\subfigure [Evolution of $w(x,t)$]
{
\includegraphics[width=6cm]{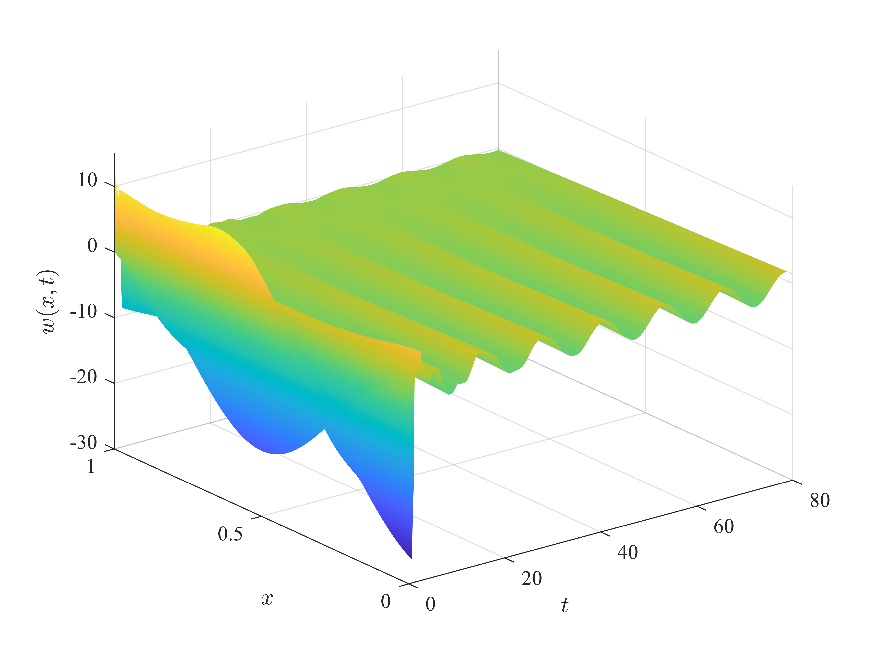}
\label{fig3}
} 
\quad 
\subfigure [Evolution of $u(t)$]
{
\includegraphics[width=6cm]{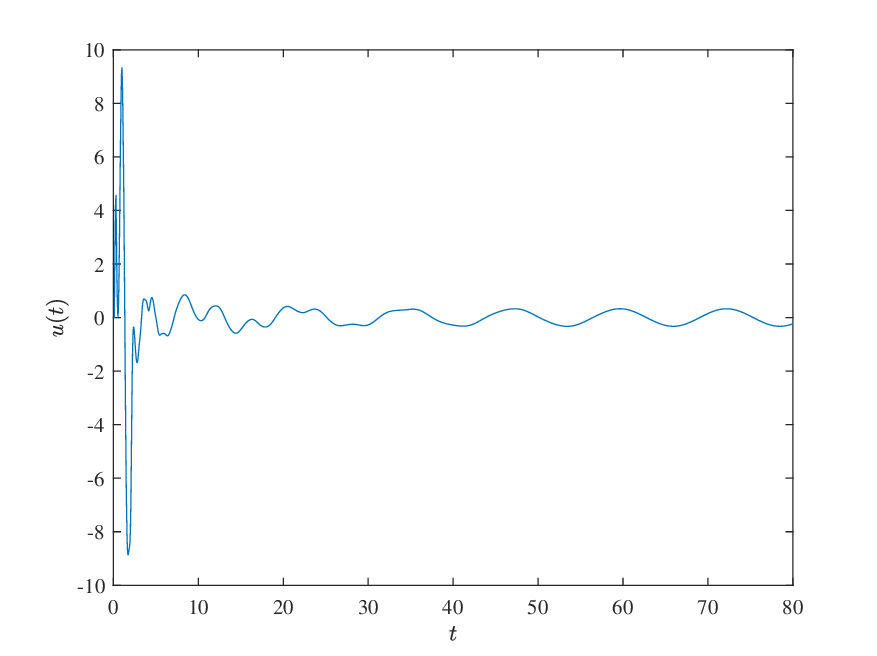}
\label{fig4}
}
\caption{Tracking performance, frequency estimate and evolution of $w(x,t)$ and $u(t)$ for $\tau=0.1$}
\label{fig5}
\end{figure}

\begin{figure}[htbp]
\centering
\subfigure[ Tracking performance $e(t)$ ]
{
\includegraphics[width=6cm]{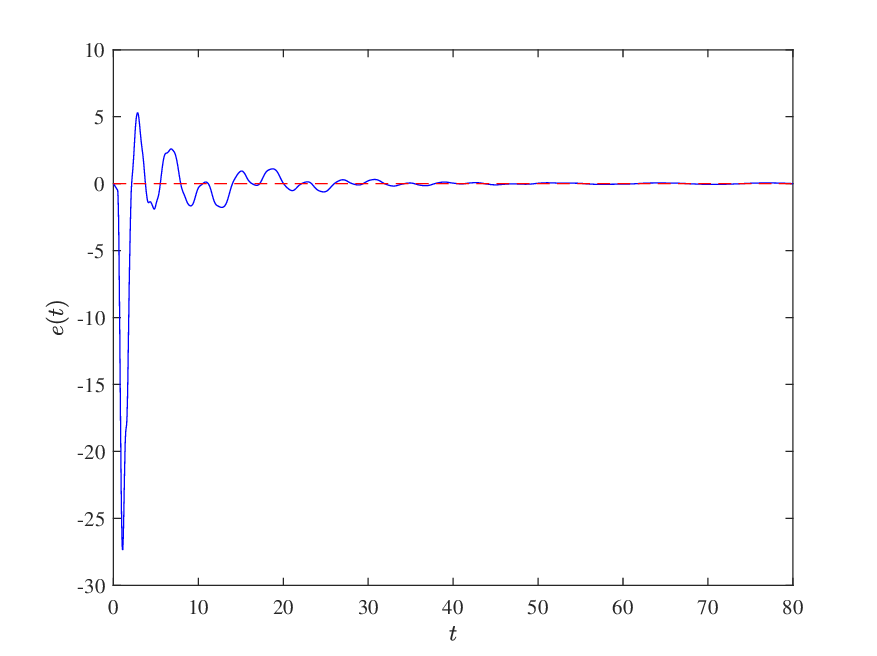}
\label{fig6}
} 
\quad 
\subfigure[ Estimate of $\theta$ ]
{
\includegraphics[width=6cm]{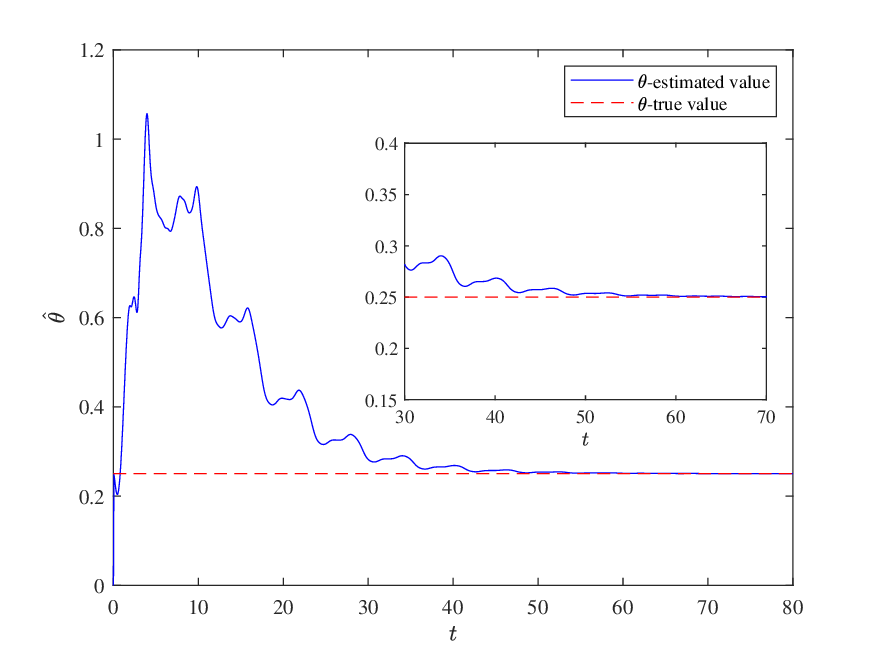}
\label{fig7}
}
\quad
\subfigure [Evolution of $w(x,t)$]
{
\includegraphics[width=6cm]{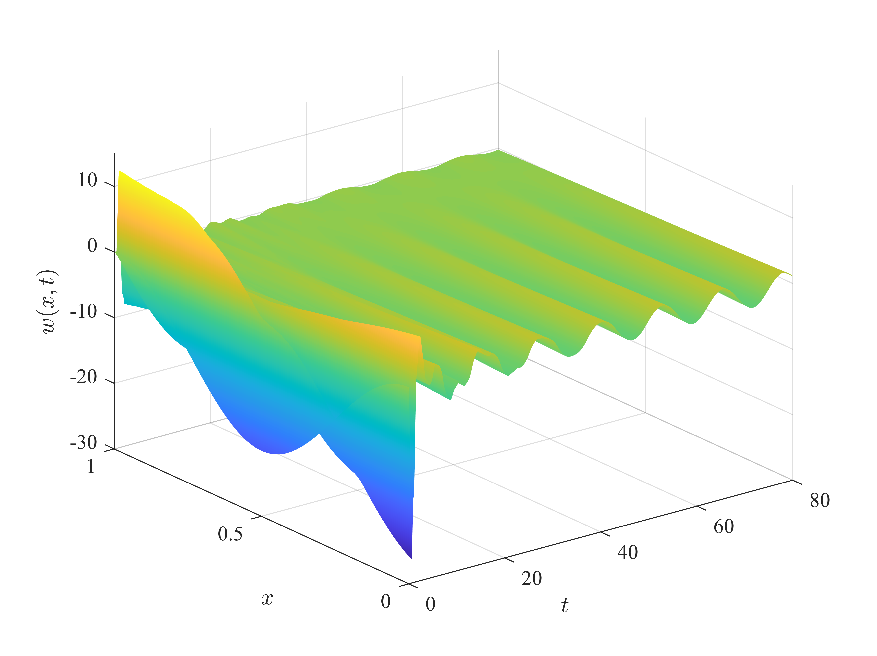}
\label{fig8}
} 
\quad 
\subfigure [Evolution of $u(t)$]
{
\includegraphics[width=6cm]{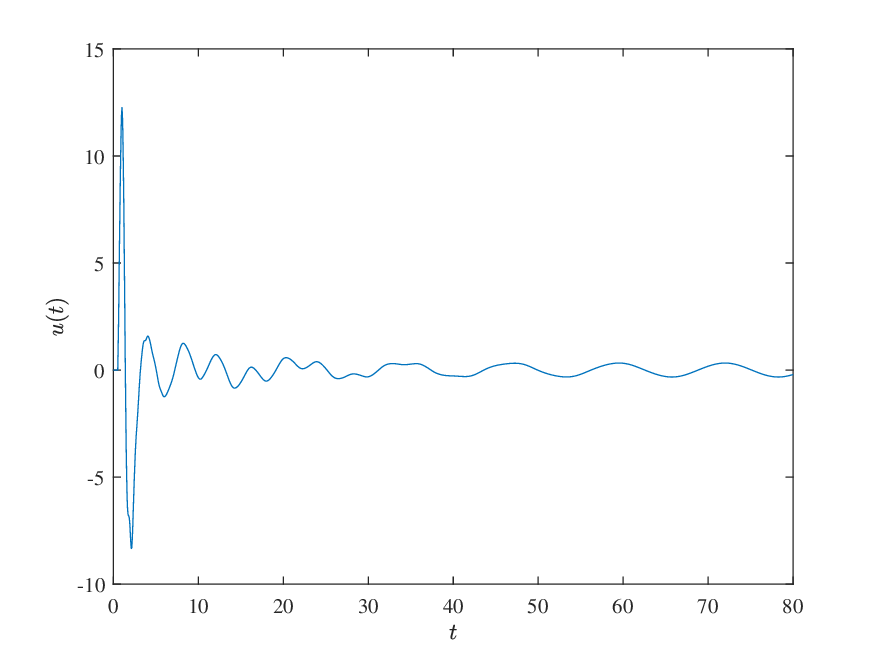}
\label{fig9}
}
\caption{Tracking performance, frequency estimate and evolution of $w(x,t)$ and $u(t)$ for $\tau=0.5$}
\label{fig10}
\end{figure}

\begin{figure}[htbp]
\centering
\subfigure[ Tracking performance $e(t)$ ]
{
\includegraphics[width=6cm]{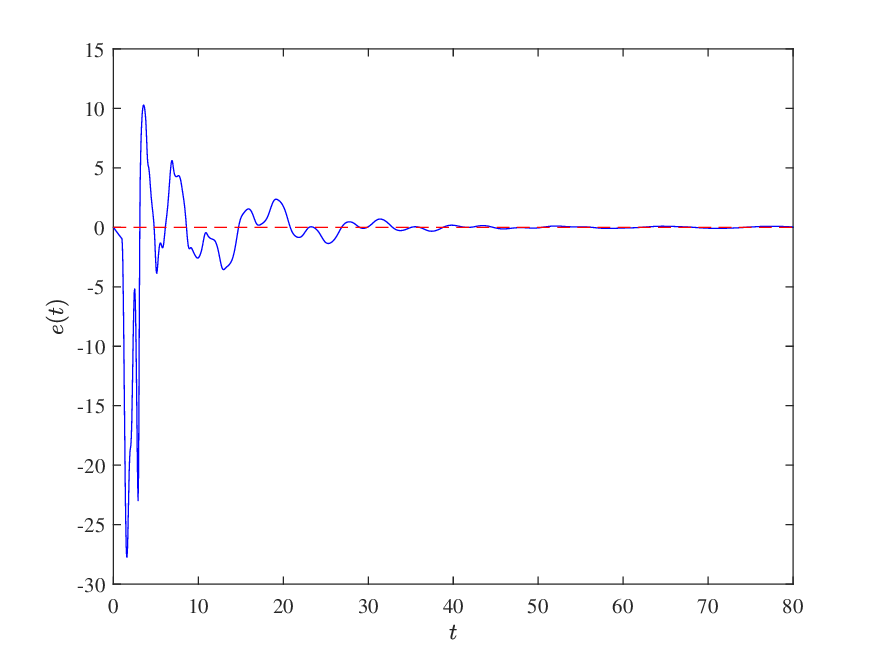}
\label{fig11}
} 
\quad 
\subfigure[ Estimate of $\theta$ ]
{
\includegraphics[width=6cm]{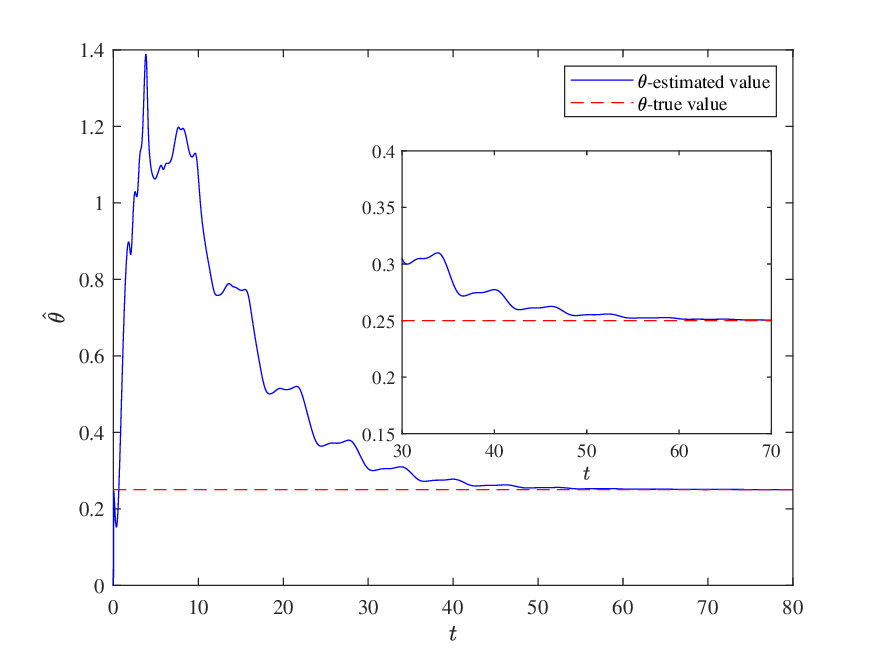}
\label{fig12}
}
\quad
\subfigure [Evolution of $w(x,t)$]
{
\includegraphics[width=6cm]{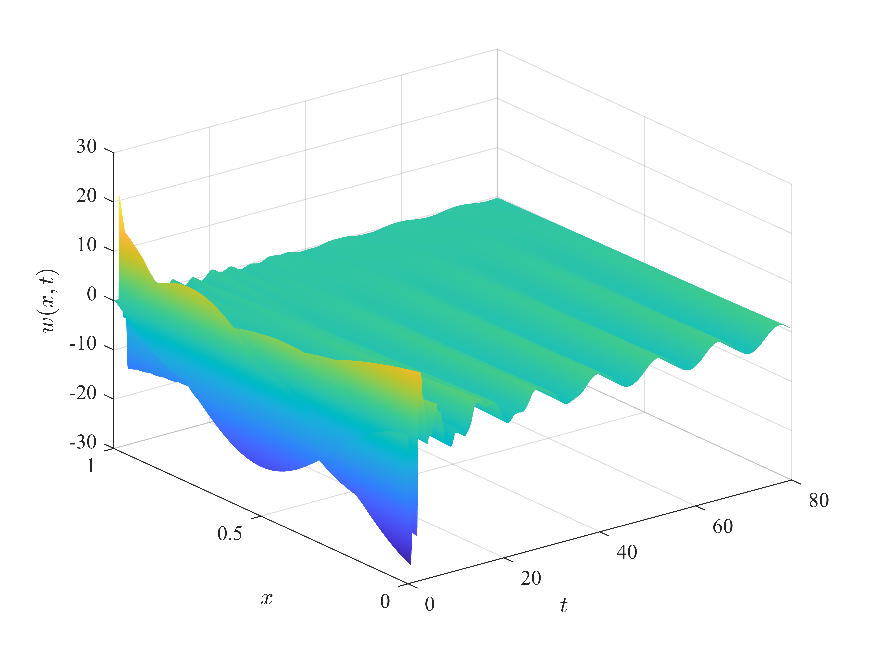}
\label{fig13}
} 
\quad 
\subfigure [Evolution of $u(t)$]
{
\includegraphics[width=6cm]{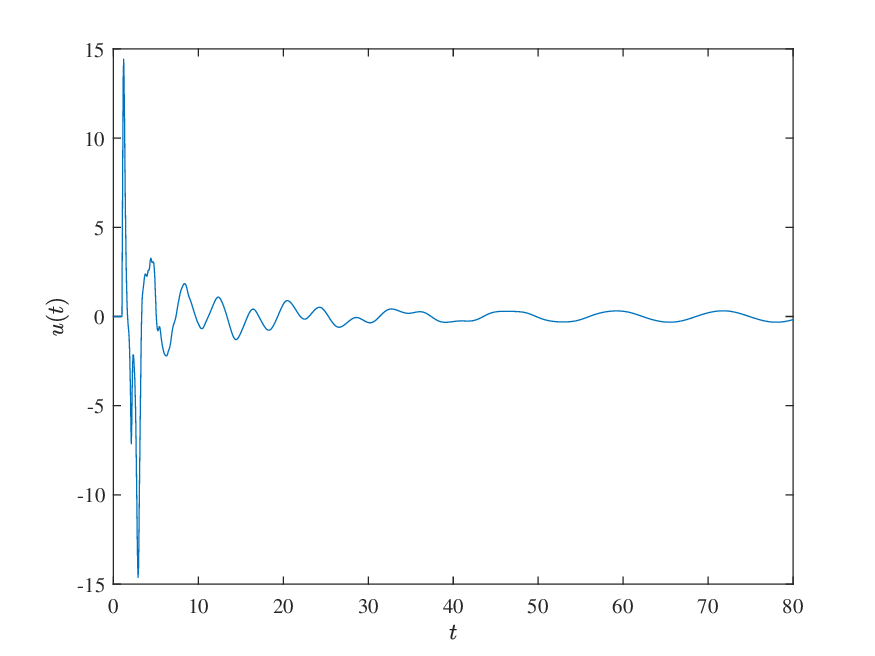}
\label{fig14}
}
\caption{Tracking performance, frequency estimate and evolution of $w(x,t)$ and $u(t)$ for $\tau=1$}
\label{fig15}
\end{figure}

In this section, we show some numerical simulations to demonstrate the effectiveness of the developed output feedback controller \eqref{69}.
Set 
$q=1, p_1(x)=2x, p_2=0, p_3=0, p_4=2, 
S=
\left[
\begin{smallmatrix}
0  & 0.25 \\
-1 & 0
\end{smallmatrix}
\right],$
and the system \eqref{1} is described as 
\begin{equation} \label{92}
\left
\{
\begin{array}{ll}
w_{t t}(x, t)=w_{x x}(x, t)+2x \sin(0.5 t),  \\
w_x(0, t)=-w(0, t), \\
w(1, t)=u(t),   \\
e(t)= y_p(t)- y_{\rm{ref}}(t)=w(0, t-\tau)-2\sin(0.5 t).
\end{array}
\right.
\end{equation}
The parameters of the controller \eqref{69} are chosed as 
$
c_0=200, c_1=1, c_2=0.1, k_0=5, k_1=10, \iota=1.
$
The initial values of the final closed-loop system \eqref{70} are taken as 
\begin{equation}
\left\{
\begin{array}{l}
w_0(x)=\hat{z}(x,0)=10(\cos(2\pi x)-1) = 0, \; 
\hat{Y}_2(x,0) = -c_2, \\
w_1(x)=Y_{10}(x)=\hat{z}_s(x,0)=\xi(0)=\hat{\chi}(0)=\hat{\varphi}(0)=\hat{\theta}(0)=0.
\end{array}
\right.
\end{equation}
The finite difference method is applied to numerically illustrate the dynamical behavior of system \eqref{70}.
Figure \ref{fig5} describes the performance of \eqref{70} when the output delay is $\tau=0.1$. 
Therein, Figure \ref{fig1} displays the convergence performance of the tracking error $e(t)$,
which demonstrate that $y_p(t)$ tracks $y_{\rm{ref}}(t)$ well after 30s.
Figure \ref{fig2} indicates the performance of the frequency estimator $\hat{\theta}$ is satisfactory.
Figure \ref{fig3} proves the boundedness of the unstable wave equation \eqref{92}, 
and Figure \ref{fig4} illustrates the development of the controller $u(t)$. 
To further demonstrate the generalizability of the proposed control strategy across different time delays, we specifically examine larger delay values ($\tau=0.5$ and $\tau=1$) compared with the case of $\tau=0.1$. 
The corresponding simulation results presented in Figures \ref{fig10} and \ref{fig15} respectively reveal that the proposed control architecture exhibits universal adaptability to arbitrarily bounded time-delay configurations.

\section{Conclusions}
This work resolves the output regulation problem for an unstable wave equation with output delay and disturbances of unknown frequency—a previously unaddressed challenge in DPS control. 
To cancel the negative effects of the output delay, we divide the time interval and introduce observers and predictors on different intervals to exponentially track the unknown states.
Regarding the unknown frequency, an adaptive observer is developed based on the adaptive internal model method to exponentially track this frequency. 
Crucially, our solution attains exponential tracking error convergence using only one measurement signal, reducing the sensing requirements of previous results. 
Moreover, this is the first result in addressing the output regulation problem for DPS with output delay, thereby expanding the scope of research on output regulation.

The principal subject for future research is to define the exosystem \eqref{2} in $\mathbb{R}^n, n\in\mathbb{N}$ to achieve a more general result. 
Another area of ongoing investigation involves applying the presented findings to other types of DPS, such as hyperbolic systems in multidimensional spatial domains and network topological structures.

\section*{CRediT authorship contribution statement}
\textbf{Shen Wang}: Writing – original draft, Investigation, Validation, Software, Methodology.
\textbf{Zhong-Jie Han}: Writing – review and editing, Supervision, Software, Methodology, Conceptualization.
\textbf{Shuangxi Huang}: Writing – review and editing, Software, Methodology.
\textbf{Zhi-Xue Zhao}: Software, Methodology.
 
\section*{Declaration of competing interest}
The authors declare that they have no known competing financial interests or personal relationships that could have appeared to influence the work reported in this paper.

\section*{Acknowledgment}
This work was supported by the Natural Science Foundation of China grant NSFC-62473281, 12326327, 12326342. 

\appendix

\section{Proof of Lemma \ref{le1}}
The system \eqref{6} can be written as
\begin{equation}\label{7}
\dfrac{d}{dx}
[\Pi (x)\quad \Pi^\prime(x)]
=[\Pi (x)\quad \Pi^\prime(x)]
\bar{S}
-p_1(x)[0 \quad E],
\end{equation} 
where
$
\bar{S}
=\left[
\begin{array}{cc}
0 & S^2 \\
E & 0
\end{array}
\right],
0,E \in \mathbb{R}^{2 \times 2}$, and $E$ is the identity matrix.
Solving \eqref{7}, we obtain
\begin{equation}\label{8}
[\Pi(x) \quad \Pi^\prime(x)]
=[p_4 e^{\tau S} \quad p_2-q p_4 e^{\tau S}] e^{\bar{S}x}
-\int_{0}^{x} [0_{1\times 2}\quad p_1(h)] e^{\bar{S}(x-h)} dh.
\end{equation}

\section{Proof of Lemma \ref{le3}}\label{sec 3}
First, according to \eqref{22}$_4$ and \eqref{22}$_5$, we have 
\begin{equation}\label{23}
g_2(x) = -c_2 g_1(0) e^{-S_\eta x}. 
\end{equation}
Regarding $S_\eta$, 
we set $\Psi_1=(1,i)^{\top}$ and $\Psi_2=(1,-i)^{\top}$ as the corresponding eigenvectors to the eigenvalues $\lambda_1 = i\omega$ and $\lambda_ 2= -i\omega$, respectively.  
Together with \eqref{23}, we right multiply $\Psi_j(x),\; j=1,2,$ to the both sides of \eqref{22}$_1$-\eqref{22}$_3$ and get
\begin{equation}\label{24}
\left\{
\begin{array}{l}
\bar{g}_{1j}^{\prime\prime}(x)
=\lambda_{j}^2\bar{g}_{1j}(x),\\
\bar{g}_{1j}^{\prime}(0)=0,\\
\bar{g}_{1j}(1)
=c_2 e^{-\lambda_j}\bar{g}_{1j}(0) - 1,
\end{array}
\right.
\end{equation}
where $\bar{g}_{1j}(x)=g_1(x)\Psi_j(x),\; j=1,2.$
We solve the solution $\bar{g}_{1j}(x)$ to \eqref{24} as
\begin{equation}\label{25}
\bar{g}_{1j}(x) =C_{1j}e^{\lambda_jx} + C_{2j}e^{-\lambda_jx},
\end{equation}
where 
\begin{equation}\label{26}
C_{1j}=C_{2j}
=\dfrac {1} {(2c_2-1)e^{-\lambda_j}-e^{\lambda_j}}
\ne 0, \; j=1,2.
\end{equation}
Finally, we obtain $g_1(x) = [ \bar{g}_{11}(x) \; \bar{g}_{12}(x) ] [ \Psi_1 \; \Psi_2 ]^{-1}$ and $g_2(x)$ by \eqref{23}.

\section{Proof of Lemma \ref{le5}}
We can find that $(S_\eta,g_1(0))$ is observable if and only if $(J,\bar{g}_1(0))$ is observable, where
$J=\Psi^{-1} S_\eta \Psi={\rm{diag}}\{i\omega,-i\omega\},$ $\bar{g}_1(0)=g_{1}(0)\Psi=(\bar{g}_{11}(0)\; \bar{g}_{12}(0))$, $\Psi=[\Psi_1 \; \Psi_2]$, 
and $\bar{g}_{11}(x), \bar{g}_{12}(x), \Psi_1, \Psi_2$ are defined in subsection \ref{sec 3}. 
By the Hautus lemma, we only need $\bar{g}_{11}(0) \ne 0$ and $\bar{g}_{12}(0)\ne0$, which is straight to see for $\omega\in(0,\infty)$ by \eqref{25} and \eqref{26}.

\section*{Data availability}
Data will be made available on request.

\bibliographystyle{elsarticle-num} 
\bibliography{reference}

\end{document}